\documentclass[]{amsart}
\usepackage{graphicx,amsfonts,amsmath,amsthm,amssymb,amscd,longtable}

\def\RP{\mathcal{RP}}
\def\P{\mathcal{P}}
\def\Q{\mathcal{Q}}
\def\U{\mathcal{U}}
\def\D{\mathcal{DRP}}
\def\tr{{\rm tr}}
\def\PSL{{\rm PSL}(2,{\mathbb C})}

\newtheorem{theorem}{Theorem}[section]
\newtheorem{lemma}[theorem]{Lemma}

\newtheorem{remark}[theorem]{Remark}

\begin{document}

\author{Elena Klimenko}
\address{Gettysburg College, Mathematics Department,
300 N. Washington St., CB 402, Gettysburg, PA 17325, USA}
\email{yklimenk@gettysburg.edu}
\author{Natalia Kopteva}
\address{LATP, UMR CNRS 6632, CMI, 39 rue F. Joliot Curie, 13453
Marseille cedex 13, FRANCE}
\email{kopteva@cmi.univ-mrs.fr}

\title[The parameter space of two-generator Kleinian groups]
{A two-dimensional slice through the parameter space of 
two-generator Kleinian groups}

\begin{abstract}
We describe all real points of the parameter
space of two-generator Kleinian groups with a parabolic generator,
that is, we describe a certain two-dimensional slice through 
this space.
In order to do this we 
gather together known discreteness criteria for two-generator
groups and present them in the form of conditions on 
parameters.
We complete the description by giving discreteness criteria 
for groups
generated by a parabolic and a $\pi$-loxodromic elements whose
commutator has real trace and present all orbifolds
uniformized by such groups.
\end{abstract}



\thanks{The first author was supported by Gettysburg College
Research and Professional Development Grant, 2005--2006.
The research of the second author was supported by 
FP6 Marie Curie IIF Fellowship and carried out at LATP (UMR CNRS 6632)}

\date{\today}

\maketitle

\section{Introduction}

A two-generator subgroup $\Gamma=\langle f,g\rangle$ of $\PSL$
is determined up to conjugacy
by its parameters
$\beta=\beta(f)={\rm tr}^2f-4$, $\beta'=\beta(g)={\rm tr}^2g-4$, and
$\gamma=\gamma(f,g)={\rm tr}[f,g]-2$ whenever $\gamma\not=0$
\cite{GM89}. So the conjugacy class of an ordered pair
$\{f,g\}$ can be identified with a point in the parameter
space ${\mathbb C}^3=\{(\beta,\beta',\gamma)\}$ whenever $\gamma\not=0$.
The subspace $\mathcal K$ of ${\mathbb C}^3$ that corresponds to 
the discrete non-elementary groups $\Gamma=\langle f,g\rangle$
is called the {\it parameter space of two-generator Kleinian
groups}.
Note that a two-generator Kleinian group $\Gamma$ 
can be represented by several points 
in $\mathcal K$, since the same group can have different generating pairs.

Among all two-generator subgroups of $\PSL$, we distinguish 
the class of $\mathcal{RP}$ {\it groups}
(two-generator groups with real parameters):
$$
\mathcal{RP}=\lbrace\Gamma :\Gamma=\langle f,g\rangle
{\rm \ for\ some\ } f,g\in{\rm PSL}(2,{\mathbb C})
{\rm \ with\ }
(\beta,\beta',\gamma)\in {\mathbb R}^3\rbrace.
$$
The aim of this paper is to completely determine 
all points in ${\mathbb C}^3$ that are parameters for
the discrete non-elementary $\RP$ groups with one generator parabolic:
$$
S_\infty=\{(\gamma,\beta):(\beta,0,\gamma) 
{\rm \ are\ parameters\ for\ some\ }
\langle f,g\rangle\in\D
\},
$$
where $\D$ denotes the class of all discrete non-elementary $\RP$ groups.
Geometrically, $S_\infty$ is a two-dimensional slice through the
six-dimensional parameter space~$\mathcal K$.

The slice $S_\infty$ intersects the well-known Riley
slice $(0,0,\gamma)$, $\gamma\in{\mathbb C}$, which consists
of all Kleinian groups generated by two parabolics.

Consider the sequence of slices $\{S_n\}_{n=2}^\infty$, where
$$
S_n=\{(\gamma,\beta):(\beta,-4\sin^2(\pi/n),\gamma) 
{\rm \ are\ parameters\ for\ some\ }
\langle f,g\rangle\in\D
\}.
$$
The first slice $S_2$ of this sequence 
is of great interest in the theory of discrete groups.
This slice consists of all
parameters for discrete $\RP$ groups with an elliptic
generator of order~2 and was investigated in \cite{GGM01}.
It was shown that if $\langle f,g\rangle$ has parameters 
$(\beta,\beta',\gamma)$, then there exists a group
$\langle f,h\rangle$ with parameters $(\beta,-4,\gamma)$
such that if $\gamma\not=0,\beta$, then $\langle f,h\rangle$
is discrete whenever $\langle f,g\rangle$ is.
Hence, the slice $S_2$ gives 
necessary discreteness conditions for a group with parameters 
$(\beta,\beta',\gamma)$, where $\beta$ and $\gamma$ are real.
It follows that every $S_n$ with $n>2$, including $S_\infty$, 
is a subset of $S_2$. 

Since a parabolic element can be viewed as the limit of a sequence
of primitive elliptic elements of order $n$ as $n\to\infty$,
the following two questions for $\{S_n\}$ and $S_\infty$
naturally arise.
\begin{itemize}
\item[(1)] Is it true that for every point $x\in S_\infty$
there exists a sequence $\{x_k\}_{k=2}^\infty$ with $x_k\in S_k$
that converges to~$x$?
\item[(2)] Is it true that for each $\varepsilon>0$ there exists
$N\in{\mathbb N}$ such that the $\varepsilon$-neighbourhood
of $S_\infty$ contains $S_n$ for all $n>N$?
\end{itemize}
Note that the structure of $S_n$ for $n>2$ is unknown.

We work out $S_\infty$ by splitting the plane $(\gamma,\beta)$
into several parts. 
It turns out that $\Gamma=\langle f,g\rangle$ 
has an invariant plane in one
of the following cases:
(1) $\gamma<0$ and $\beta\leq -4$;
(2) $\gamma>0$ and $\beta\geq -4$. 
Such discrete groups 
were investigated, for example, in \cite{KS98} and
\cite{GiM91,Kna68,Mat82}, respectively.
If $\gamma<0$ and $\beta>-4$, then $\Gamma$ is truly spatial
(non-elementary and without invariant plane)
and this case is treated in~\cite{KK05}.
We get these dicreteness criteria together
and transform them
into conditions on $\beta$ and $\gamma$ if it was not done before.

So the last case to consider is when $\gamma>0$ and $\beta<-4$.
In this case $\Gamma$ is truly spatial with $f$ $\pi$-loxodromic.
We complete the study of the slice 
$S_\infty$ by giving discreteness criteria for 
such groups.

The paper is organised as follows.
In Section~2, discreteness criteria are given for 
truly spatial $\RP$ groups $\Gamma$ 
generated by a $\pi$-loxodromic and a parabolic elements
(Theorems~\ref{criterion_psl} and~\ref{criterion_par}).
In Section~3,
for each such discrete $\Gamma$ we obtain a presentation and the
Kleinian orbifold $Q(\Gamma)$ (Theorem~\ref{groups}).
Section~4 is devoted to the analysis of the parameter space.
We completely describe the slice $S_\infty$ by 
giving explicit formulas for the parameters $\beta$ and $\gamma$.
We also program the obtained formulas in the package Maple 7.0
and
plot a part of $S_\infty$ on the $(\gamma,\beta)$-plane to give an idea
of how it looks like.

\section{Discreteness criteria}

Recall that an element $f\in\PSL$ with real $\beta(f)$
is
{\it elliptic}, {\it parabolic}, {\it hyperbolic}, or
{\it $\pi$-loxodromic}
according to whether
$\beta(f)\in[-4,0)$, $\beta(f)=0$, $\beta(f)\in(0,+\infty)$, or
$\beta(f)\in(-\infty,-4)$.
If $\beta(f)\notin[-4,+\infty)$,
then $f$ is called {\it strictly loxodromic}.

An elliptic element $f$ of order $n$ is said to be {\it non-primitive}
if it is a rotation through $2\pi q/n$, where $q$ and $n$ are
coprime ($1<q<n/2$). If $f$ is a rotation through $2\pi/n$,
then it is called {\it primitive}.

\begin{theorem}\label{criterion_psl}
Let $f\in\PSL$ be a $\pi$-loxodromic element,
$g\in\PSL$ be a parabolic element,
and let $\Gamma=\langle f,g\rangle$ be a non-elementary
$\RP$ group without invariant plane.
Then
\begin{itemize}
\item[(1)]
there  exist unique elements $h_1,h_2\in\PSL$ such that
$h_1^2=fg^{-1}f^{-1}g^{-1}$ and $(h_1g)^2=1$,
$h_2^2=f^{-1}g^{-1}f^2gf^{-1}$ and $(h_2fg^{-1}f^{-1})^2=1$.
\item[(2)]
the group $\Gamma$ is discrete if and only if one of the following
conditions holds:
\begin{itemize}
\item[(i)]
$h_1$ is either a hyperbolic, or parabolic, or primitive elliptic
element of even order $m\geq 4$, and
$h_2$ is either a hyperbolic, or parabolic, or primitive elliptic
element of order $p\geq 3$;
\item[(ii)]
$h_1$ is a primitive elliptic
element of odd order $m\geq 3$, and $h_2h_1$ is either a hyperbolic,
or parabolic, or primitive elliptic element of order $k\geq 3$.
\end{itemize}
\end{itemize}
\end{theorem}

\subsection*{Basic geometric construction}

We will construct a group $\Gamma^*$ that
contains $\Gamma=\langle f,g\rangle$ as a subgroup of finite index.
The idea is to find $\Gamma^*$ so that a fundamental polyhedron for
a discrete $\Gamma^*$ can be easily constructed.
It will be clear from the construction that $\Gamma$ is
commensurable with a reflection group which either coincides
with $\Gamma^*$ or is an index 2 subgroup of $\Gamma^*$.
The construction
presented below will be used throughout Sections~2 and~3
and we shall use the notation introduced here.

\bigskip

Let $f$ and $g$ be as in the statement of Theorem~\ref{criterion_psl}.
Since $\Gamma$ is a non-elementary $\RP$ group without invariant plane,
there exists an invariant plane of $g$, say $\eta$,
which is orthogonal to the axis of $f$ \cite[Theorem~2]{KK02}.

Denote by $M$ the fixed point of $g$ and by $\omega$ the plane
that passes through $M$ and $f$ (we denote elements and their
axes by the same letters when it does not lead to any confusion).
Note that $f$ keeps $\omega$ invariant.
Since $f$ is orthogonal to $\eta$, $\omega$ is also orthogonal
to $\eta$. Let $e$ be the half-turn with the axis $\omega\cap\eta$.
Then $e$ passes through $M$ and is orthogonal to~$f$.

Let $e_f$ and $e_g$ be half-turns such that
\begin{equation}\label{efeg}
f=e_fe {\rm \quad and \quad} g=e_ge.
\end{equation}
Then $e_f$ is orthogonal to $\omega$ and $e_g$ lies in $\eta$.

Let $\tau$ be the plane passing through $e_g$ orthogonally to $\eta$
and let $\sigma=e_f(\tau)$.
The planes $\tau$ and $\omega$ are parallel and $M$ is their common
point on the boundary $\partial{\mathbb H}^3$.
Since $e_f$ is orthogonal to $\omega$, the planes
$\sigma$ and $\omega$ are also parallel
with the common point $e_f(M)$ on $\partial{\mathbb H}^3$.
Since $e_f(M)\not=M$, the planes $\omega$, $\sigma$, and $\tau$
do not have a common point in
$\overline{{\mathbb H}^3}={\mathbb H}^3\cup\partial{\mathbb H}^3$.
Therefore,
there exists a unique plane $\delta$ orthogonal to all $\omega$,
$\sigma$, and $\tau$.
It is clear that $e_f\subset\delta$.

Consider two extensions of $\Gamma$:
$\widetilde\Gamma=\langle f,g,e\rangle$ and
$\Gamma^*=\langle f,g,e,R_\omega\rangle$.
(We denote the reflection in a plane $\kappa$ by $R_\kappa$.)
One can show that $\widetilde\Gamma=\langle e_f,e_g,e\rangle$
and $\Gamma^*=\langle e_f,R_\eta,R_\omega,R_\tau\rangle$.
From (\ref{efeg}), it follows that $\widetilde\Gamma$ contains $\Gamma$
as a subgroup of index at most~2.
Moreover, $\widetilde\Gamma$ is the orientation preserving subgroup of
$\Gamma^*$ and, hence, $\Gamma^*$ contains $\Gamma$ as a subgroup
of finite index. Therefore, $\Gamma$, $\widetilde\Gamma$, and
$\Gamma^*$ are either all discrete, or all non-discrete.
We then concentrate on the group~$\Gamma^*$.

\begin{figure}[htbp]
\centering
\includegraphics[width=10 cm]{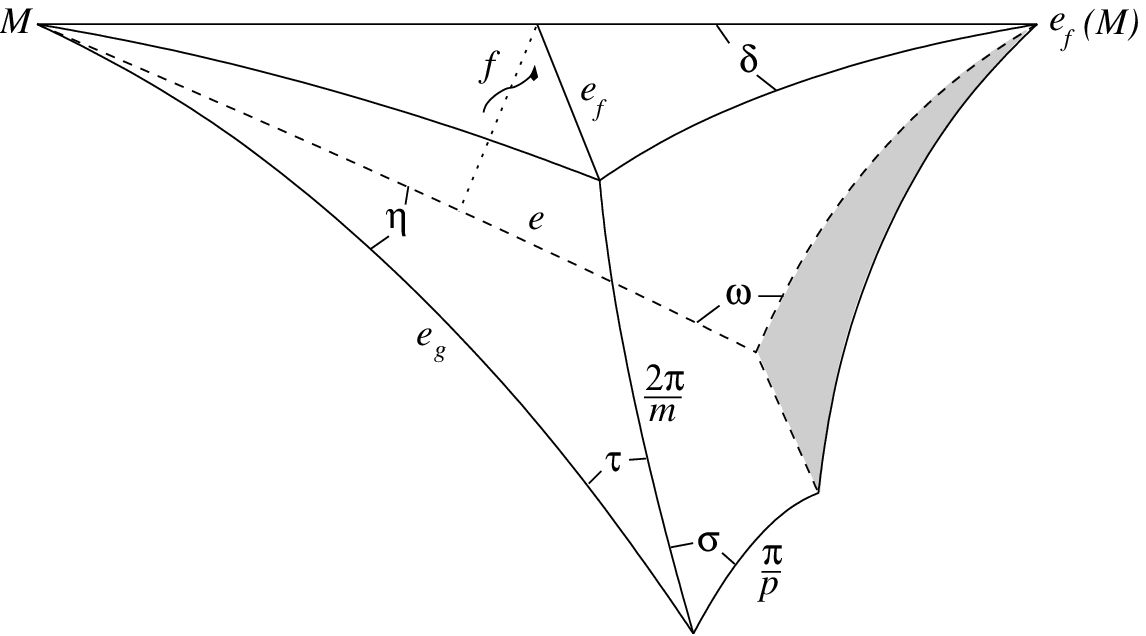}
\caption{Polyhedron $\P^*$}\label{fund_poly1}
\end{figure}

Let $\P^*$ be the infinite volume polyhedron bounded by
$\eta$, $\omega$, $\tau$, $\sigma$, and $\delta$. $\P^*$ has five right
dihedral angles (between faces lying in $\eta$ and $\omega$,
$\eta$ and $\tau$, $\delta$ and $\omega$, $\delta$ and $\tau$,
and $\delta$ and $\sigma$).
The plane $\sigma$ may either intersect with, or be parallel to, or
be disjoint from each of $\tau$ and $\eta$.

If $\sigma$ and $\tau$ intersect, then we denote the dihedral angle of
$\P^*$ between them by $2\pi/m$, where $m>2$, $m$ is not necessary an
integer. We keep the notation $2\pi/m$ taking $m=\infty$ and
$m=\overline\infty$ for parallel or disjoint $\sigma$ and $\tau$,
respectively.
Similarly, we denote the ``dihedral angle'' between
$\eta$ and $\sigma$ by $\pi/p$, where $p>2$ is real,
$\infty$, or $\overline\infty$.
(We regard $\overline\infty>\infty>x$, $x/\infty=x/\overline\infty=0$,
$\infty/x=\infty$, $\overline\infty/x=\overline\infty$ for any
positive real $x$.)
$\P^*$ exists in ${\mathbb H}^3$ for all $m>2$ and $p>2$
by~\cite{Vin85}.

In Figure~\ref{fund_poly1}, $\P^*$ is drawn under assumption that
$m<\infty$, $p<\infty$, and
$1/2+1/p+2/m>1$. The shaded triangle shows the hyperbolic plane
orthogonal to $\eta$, $\sigma$, and $\omega$. Note that
this plane is not a face of $\P^*$ and is shown only to underline
the combinatorial structure of $\P^*$. In figures, we do not label
dihedral angles of $\pi/2$ in order to not overload the picture.

Suppose now that $m<\infty$, that is $\sigma$ and $\tau$ intersect.
Let $\xi$ be the plane passing through $e_f$ orthogonally to $\delta$.
Then $\xi$ is orthogonal to $\omega$. One can see that
$\sigma=R_\xi(\tau)$ and
$\xi$ is the bisector of the dihedral angle of $\P^*$ made by
$\tau$ and~$\sigma$.

\begin{figure}[htbp]
\centering
\includegraphics[width=6 cm]{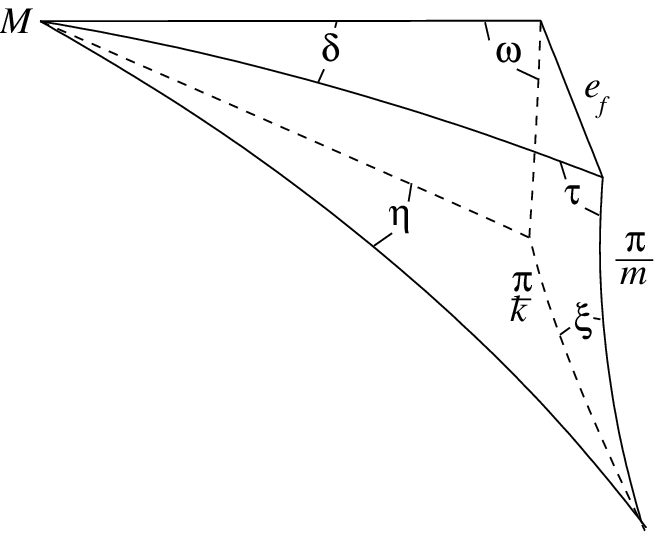}
\caption{Polyhedron $\Q^*$}\label{fund_poly2}
\end{figure}

Let $\Q^*$ be the polyhedron bounded by $\eta$, $\tau$, $\omega$,
$\delta$, and $\xi$. $\Q^*$ has six dihedral angles of $\pi/2$;
the dihedral angle between $\tau$ and $\xi$ is equal to $\pi/m$ with
$2<m<\infty$.
Denote the ``dihedral angle'' between $\eta$ and $\xi$ by $\pi/k$,
where $k>2$ is real, $k=\infty$, or $k=\overline\infty$.
$\Q^*$ exists in ${\mathbb H}^3$ for all $m>2$ and $k>2$
by~\cite{Vin85}.
Note that $R_\xi$ is not necessary in $\Gamma^*$, but if it is and
if $\Gamma^*$ is discrete, then we will see that $\Q^*$ is a
fundamental polyhedron for~$\Gamma^*$.
In Figure~\ref{fund_poly2}, $\Q^*$ is drawn under assumption that
$1/2+1/k+1/m>1$.

\begin{lemma}\label{hs}
Let $f\in\PSL$ be a $\pi$-loxodromic element,
$g\in\PSL$ be a parabolic element,
and let $\Gamma=\langle f,g\rangle$ be a non-elementary
$\RP$ group without invariant plane.
Then there  exist unique elements
$h_1,h_2\in\PSL$ such that
\begin{itemize}
\item[(1)] $h_1^2=fg^{-1}f^{-1}g^{-1}$ and $(h_1g)^2=1$,
\item[(2)] $h_2^2=f^{-1}g^{-1}f^2gf^{-1}$ and $(h_2fg^{-1}f^{-1})^2=1$.
\end{itemize}
Moreover, the elements $h_1$ and $h_2$ are not strictly loxodromic.
\end{lemma}

\begin{proof}
First, note that $R_\sigma=e_fR_\tau e_f$
and $g=R_\tau R_\omega$.
Therefore,
\begin{equation}\label{fprim}
R_\sigma R_\omega= e_f R_\tau e_f R_\omega=e_fR_\tau R_\omega e_f=
e_fge_f=fg^{-1}f^{-1}.
\end{equation}

Let us show that if we take $h_1=R_\xi R_\tau=R_\sigma R_\xi$, then
the assertion (1) of the lemma hold. Indeed,
$$
h_1^2=R_\sigma R_\tau=(R_\sigma R_\omega)(R_\omega R_\tau)=fg^{-1}f^{-1}g^{-1}.
$$
Moreover, $h_1g=(R_\xi R_\tau)(R_\tau R_\omega)=R_\xi R_\omega$. Since
$\xi$ and $\omega$ are orthogonal, $(R_\xi R_\omega)^2=1$.
Hence, $(h_1g)^2=1$.
Note also that since $h_1$ is a product of two reflections,
$h_1$ is not strictly loxodromic.

Now let us show that $h_1$ is unique.
The element $fg^{-1}f^{-1}g^{-1}$ is uniquely determined
as an element of $\PSL$.

If $fg^{-1}f^{-1}g^{-1}$ is parabolic, it has only one square
root~$h_1$. Suppose that $fg^{-1}f^{-1}g^{-1}$ is hyperbolic. Then
it has exactly two square roots, one of which is $h_1$ defined above
and the other, denoted $\overline h_1$, is a $\pi$-loxodromic element
with the same axis and translation length as $h_1$. Clearly,
$(\overline h_1g)^2\not=1$.

If $fg^{-1}f^{-1}g^{-1}$ is elliptic, then it also has two square roots
$h_1$ and $\overline h_1$,
both are elliptic elements. The element $\overline h_1$ is elliptic
with the same axis as $h_1$ and with rotation angle $(\pi-2\pi/m)$,
while $h_1$ is a rotation through $2\pi/m$ in the opposite direction.
Again, $(\overline h_1g)^2\not=1$.

Now we take
$$
h_2=R_\eta R_\sigma=(R_\eta R_\tau)(R_\tau R_\sigma)=e_gh_1^{-2}=
efgf^{-1}.
$$
Then
$$
h_2^2=f^{-1}g^{-1}f^2gf^{-1} {\rm\quad and\quad} (fg^{-1}f^{-1}h_2)^2=1.
$$
These two conditions determine
$h_2$ uniquely.
\end{proof}

Note that the elements $h_1, h_2$ defined in Lemma~\ref{hs} determine
combinatorial and metric structures of $\P^*$. For example, if
$h_1$ is elliptic, then its rotation angle is equal to the dihedral angle
of $\P^*$ between $\sigma$ and $\tau$. If $h_2$ is elliptic, then its
rotation angle is equal to the doubled dihedral angle of $\P^*$ between
$\eta$ and $\sigma$. Vice versa, if the metric structure of $\P^*$
is fixed, then the types of elements $h_1$ and $h_2$ can be determined.

The same can be said about $\Q^*$ and the elements $h_1$ and $h_2h_1$.
The element $h_2h_1$ is responsible for the
mutual position of the planes $\eta$ and~$\xi$ (see the proof of 
Lemma~\ref{h1h2}).

Lemmas~\ref{h1}--\ref{h1h2} below give some necessary conditions
for discreteness of $\Gamma$ via conditions on elements $h_1$ and $h_2$. 
One needs to
keep in mind the connection between these elements and the
polyhedra $\P^*$
and~$\Q^*$.

\begin{lemma}\label{h1}
If $\Gamma$ is discrete, then
$h_1$ is either a hyperbolic, or parabolic, or primitive elliptic
element of order $m\geq 3$.
\end{lemma}

\begin{proof}
The subgroup $H=\langle g,fgf^{-1}\rangle$ of $\Gamma$
keeps $\delta$ invariant and is conjugate to a
subgroup of ${\rm PSL}(2,\mathbb R)$.
Since $\Gamma$ is discrete, $H$ must be discrete.
By \cite{Mat82} or \cite{Bea88},
the group $H$ is discrete if and only if either

(1) $fg^{-1}f^{-1}g^{-1}=h_1^2$
is a hyperbolic, or a parabolic, or a primitive elliptic
element, or 

(2) $h_1$ is a primitive elliptic element of odd order
$m$, where $m\geq 3$.

If $h_1^2$ is parabolic of hyperbolic, then $h_1$
is parabolic or hyperbolic, respectively.
If $h_1^2$ is a primitive elliptic element,
then $h_1$ is a primitive elliptic of even order $m\geq 4$.
\end{proof}

\begin{lemma}\label{h2}
If $\Gamma$ is discrete,
then $h_2$ is either a hyperbolic, or parabolic, or primitive elliptic
element of order $p\geq 3$.
\end{lemma}

\begin{proof}
Let $\kappa$ be the plane orthogonal to
$\eta$, $\sigma$, and $\omega$.
The subgroup $H=\langle e,fgf^{-1}\rangle$ of $\widetilde\Gamma$
keeps the plane $\kappa$ invariant and is conjugate to a
subgroup of ${\rm PSL}(2,\mathbb R)$.
By \cite{Mat82},
$H$ is discrete if and only if $h_2=efgf^{-1}$
is either a hyperbolic, or parabolic, or primitive elliptic
element of order $p\geq 3$.
\end{proof}

\begin{lemma}\label{h1h2}
If $\Gamma$ is discrete and $h_1$ is a primitive elliptic element
of odd order, then $h_2h_1$ is either a hyperbolic,
or parabolic, or primitive elliptic element of order $k\geq 3$.
\end{lemma}

\begin{proof}
Recall that
$\Gamma^*=\langle e_f,R_\eta,R_\tau,R_\omega\rangle$.
Since $h_1$ has odd order and $h_1^2\in\Gamma^*$,
$h_1\in\Gamma^*$. Since, moreover,
$h_1=R_\xi R_\tau$,
$e_f=R_\delta R_\xi$, and $R_\tau\in\Gamma^*$,
both $R_\xi$ and $R_\delta$ are also in~$\Gamma^*$.
Further, since the plane $\xi$ is orthogonal to $\omega$, the group
$\langle R_\eta R_\delta, e_f\rangle$
keeps $\omega$ invariant and is conjugate to a
subgroup of ${\rm PSL}(2,\mathbb R)$.
It is clear that $\langle R_\eta R_\delta, e_f\rangle$ is
discrete if and only if $R_\eta R_\xi=h_2h_1$ is a hyperbolic,
parabolic, or primitive elliptic element of order $k\geq 3$ \cite{Mat82}.
\end{proof}

\medskip
\noindent
{\it Proof of Theorem~\ref{criterion_psl}.}
Lemma~\ref{hs} proves existence and uniqueness of elements $h_1$
and $h_2$. Now we prove part (2) of the theorem.

If $\Gamma$ is discrete then
$h_1$ is either a hyperbolic,
or parabolic, or primitive elliptic element of order $m\geq 3$
by Lemma~\ref{h1}.
We split the discrete groups $\Gamma$ into two families.
The first family consists of those groups for
which $h_1$ is hyperbolic, parabolic, or primitive elliptic
of even order.
By Lemma~\ref{h2}, for these groups $h_2$ is a hyperbolic,
parabolic, or primitive elliptic element.

The second family consists of the discrete groups with
$h_1$ elliptic of odd order.
Then by Lemma~\ref{h1h2}, $h_2h_1$ is a hyperbolic, or parabolic,
or primitive elliptic element of order $k\geq 3$.
(Note that in this case 
$h_2$ is necessarily hyperbolic or primitive elliptic.)

So if $\Gamma$ is discrete, then either (2)(i) or
(2)(ii) of Theorem~\ref{criterion_psl} can occur.
Clearly, if neither (2)(i) nor (2)(ii) holds, then $\Gamma$ is
not discrete by Lemmas~\ref{h1}--\ref{h1h2}.

\smallskip
Now prove that each of (2)(i) and (2)(ii) is a sufficient condition
for $\Gamma$ to be discrete.
In each of the two cases we will give a fundamental polyhedron for
$\Gamma^*$ to show, by using the Poincar\'e polyhedron theorem
\cite{EP94}, that
$\Gamma^*$ is discrete.

Suppose that (2)(i) holds. Then since $m$ is even,
the group $G_1$ generated by the side pairing transformations
$R_\eta$, $R_\omega$,
$R_\sigma$, $R_\tau$, and $e_f$
and the polyhedron $\P^*$ satisfy the 
Poincar\'e polyhedron theorem, $G_1$ is discrete and $\P^*$
is its fundamental polyhedron.
Obviously, $G_1=\Gamma^*$.

Suppose that (2)(ii) holds. Then the
group $G_2$ generated by the side pairing transformations
$R_\eta$, $R_\omega$,
$R_\xi$, $R_\tau$, and $R_\delta$ and
the polyhedron $\Q^*$ satisfy the Poincar\'e theorem, $G_2$ is discrete,
and $\Q^*$ is its fundamental polyhedron.

In the proof of Lemma~\ref{h1h2} it was shown that, for $m$ odd,
$R_\xi\in\Gamma^*$ and $R_\delta\in\Gamma^*$. 
Moreover, $e_f=R_\xi R_\delta$.
Hence, $G_2=\Gamma^*$, so $\Gamma^*$ is discrete.

Theorem~\ref{criterion_psl} is proved.\qed

\bigskip

Our next goal is to compute parameters $(\beta(f),\beta(g),\gamma(f,g))$
for both series of discrete groups listed in Theorem~\ref{criterion_psl}.

If $f\in\PSL$ is a loxodromic element with translation length $d_f$ and
rotation angle $\theta_f$, then
$$
\tr^2f=4\cosh^2 \frac{d_f+i\theta_f}2
$$
and $\lambda_f=d_f+i\theta_f$ is called the
{\it complex translation length} of~$f$.

Note that if $f$ is hyperbolic then $\theta_f=0$ and $\tr^2f=4\cosh^2(d_f/2)$.
If $f$ is elliptic then $d_f=0$ and $\tr^2f=4\cos^2(\theta_f/2)$.
If $f$ is parabolic then $\tr^2f=4$; 
by convention we set $d_f=\theta_f=0$.

We define the set
$$
\mathcal{U}=\{u:u=i\pi/p
{\rm \ for\ some\ } p\in{\mathbb Z}, p\geq 2\}\cup[0,+\infty).
$$
In other words, the set $\U$ consists of 
all complex translation half-lengths
$u=\lambda_f/2$ for hyperbolic,
parabolic, and primitive elliptic elements~$f$.
Furthermore, we define a function
$t:\mathcal{U}\to\{2,3,4,\dots\}\cup\{\infty,\overline\infty\}$
as follows:
$$
t(u)=\left\{
\begin{array}{lll}
p & {\rm if} & u=i\pi/p,\\
\infty & {\rm if} & u=0,\\
\overline\infty & {\rm if} & u\in(0,+\infty).
\end{array}
\right.
$$

Given $u\in\U$ and $f$ with $\tr^2f=4\cosh^2u$, 
$t(u)$ determines the type of $f$ and, moreover, its
order if $f$ is elliptic.
Note also that since we regard $\infty/n=\infty$ and
$\overline\infty/n=\overline\infty$, an expression of the form
$(t(u),n)=1$ with $n>1$ means, in particular, that $t(u)$ is finite.

\begin{theorem}\label{criterion_par}
Let $f,g\in\PSL$ with $\beta(f)<-4$,
$\beta(g)=0$, and $\gamma(f,g)>0$.
Then $\Gamma=\langle f,g\rangle$
is discrete if and only if one of the following holds:
\begin{enumerate}
\item $\gamma(f,g)=4\cosh^2u$ and $\beta(f)=-4\cosh^2v/\gamma(f,g)-4$, 
where $u,v\in\U$ with $t(u)\geq 4$, $(t(u),2)=2$, and $t(v)\geq 3$;
\item $\gamma(f,g)=4\cosh^2u$ and $\beta(f)=-4\cosh^2v-4$, where
$u,v\in\U$ with $t(u)\geq 3$, $(t(u),2)=1$, and $t(v)\geq 3$.
\end{enumerate}
\end{theorem}

\begin{proof}
Obviously, $\beta(f)<-4$ and $\beta(g)=0$ if and only if
$f$ is $\pi$-loxodromic and $g$ is parabolic.
With this choice of $\beta(f)$ and $\beta(g)$, 
$\gamma(f,g)>0$ if and only if the group
$\Gamma=\langle f,g\rangle$ is a non-elementary $\RP$ group without 
invariant plane~\cite{KK02}.
This means that the hypotheses of Theorem~\ref{criterion_par}
are equivalent
to the hypotheses of Theorem~\ref{criterion_psl}.
Therefore, in order to prove Theorem~\ref{criterion_par} it is
sufficient to
calculate the parameters $\beta(f)$ and $\gamma(f,g)$ for
both families of the discrete groups listed in
Theorem~\ref{criterion_psl}.

Let $\sigma'$ be the image of $\sigma$ under $R_\omega$, that is
$R_{\sigma'}=R_\omega R_\sigma R_\omega$.
Using the identity~(\ref{fprim})
and the fact that $g=R_\tau R_\omega$,
we have
$$
[f,g]=fgf^{-1}g^{-1}=(R_\omega R_\sigma)(R_\omega R_\tau)=
(R_{\sigma'} R_\omega)(R_\omega R_\tau)=R_{\sigma'}R_\tau.
$$
Note that $\sigma'$
and $\tau$ are disjoint and $\delta$ is orthogonal to both of them.
Therefore, $[f,g]$ is a hyperbolic element with the
axis lying in $\delta$ and the translation length $2d$, where $d$ is
the distance between $\sigma'$ and $\tau$.
Hence, since $\gamma(f,g)>0$,
$$\gamma(f,g)=\tr[f,g]-2=+2\cosh d-2.$$

\begin{figure}[htbp]
\centering
\includegraphics[width=10 cm]{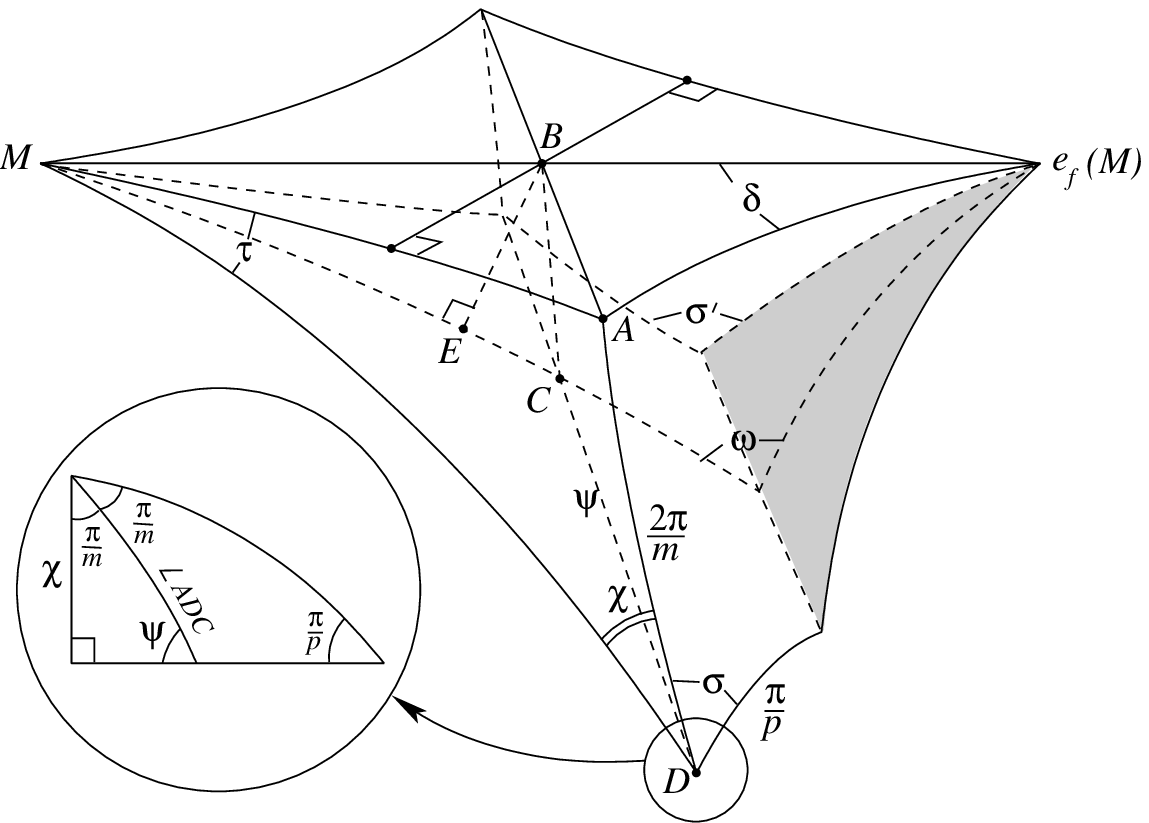}
\caption{}\label{poly_par}
\end{figure}

Now, using generalised triangles in
the plane $\delta$, it is not difficult to calculate that
$$
\gamma(f,g)=\left\{
\begin{array}{lll}
4\cos^2(\pi/m) & {\rm if} & 3\leq m<\infty,\\
4 & {\rm if} & m=\infty,\\
4\cosh^2(d(\sigma,\tau)/2) & {\rm if} & m=\overline\infty,\\
\end{array}
\right.
$$
where $d(\sigma,\tau)$ is the distance between $\sigma$ and $\tau$ if they are
disjoint.
Hence, 
$$\gamma(f,g)=4\cosh^2u,$$ 
where $u\in\U$, $t(u)=m\geq 3$.

Let us calculate $\beta(f)$.
The element $f$ is $\pi$-loxodromic if and only if
$\tr^2f=4\cosh^2(T+i\pi/2)=-4\sinh^2T$, where $2T$ is the
translation length of $f$. That is,
$$
\beta(f)=-4\sinh^2T-4.
$$
Note that $T$ is the distance between $e$ and $e_f$. It is measured in
$\omega$ and equals $BE$ (see Figure~\ref{poly_par}).

Suppose that we are in case (2)(i) of Theorem~\ref{criterion_psl},
that is $(t(u),2)=2$,
and that
$\sigma$ and $\tau$ intersect.
Recall that $\xi$ is the bisector of the dihedral angle of $\P^*$ made
by $\sigma$ and $\tau$.
Let $\psi$ be the angle that $\xi$ makes with $\eta$. Note that
$\psi=\angle BCE$.
From the link of $D$, we have that
$$
\cos\chi=\frac{\cos(\pi/p)}{\sin(2\pi/m)}=\frac{\cos\psi}{\sin(\pi/m)}
$$
and, therefore,
\begin{equation}\label{psi}
\cos\psi=\frac{\cos(\pi/p)}{2\cos(\pi/m)}.
\end{equation}
Further, from the link of $D$,
\begin{equation}\label{adc}
\cos\angle ADC=
\frac{\cos\psi\cdot\cos(\pi/m)}{\sin\psi\cdot\sin(\pi/m)}.
\end{equation}
From the $\triangle ABM$, $\cosh^2 AB=1/\sin(\pi/m)$
and, from the quadrilateral $ABCD$,
\begin{equation}\label{bc}
\sinh BC=\frac{\cos \angle ADC}{\sinh AB}
\end{equation}
Finally, from $\triangle BCE$,
\begin{equation}\label{sinht}
\sinh T=\sinh BE=\sin\psi\cdot \sinh BC.
\end{equation}
Combining (\ref{psi})--(\ref{sinht}), we have that
$$
\sinh^2T=\frac{\cos^2(\pi/p)}{4\cos^2(\pi/m)}=
\frac{\cos^2(\pi/p)}{\gamma(f,g)}.
$$
Similar calculations can be done for parallel or disjoint $\sigma$
and $\tau$.
Hence,
$\beta(f)=-\sinh^2T-4=-\cosh^2v/\gamma(f,g)-4$, 
where $v\in\U$, $t(v)\geq 3$.

Now note that in case (2)(ii) of Theorem~\ref{criterion_psl},
the angle $\psi=\angle BCE$ must be of the form $\pi/k$, $k\geq 3$
is an integer, $\infty$, or $\overline\infty$.
Then we need to recompute the formulas (\ref{adc})--(\ref{sinht})
with $\psi=\pi/k$:
$$
\cos\angle ADC=\frac{\cos(\pi/k)\cdot\cos(\pi/m)}
{\sin(\pi/k)\cdot\sin(\pi/m)},
\quad
\sinh BC=\frac{\cos\phi}{\sinh a}=\frac{\cos(\pi/k)}{\sin(\pi/k)}.
$$
Then
$$
\sinh T=\sin\psi \cdot\sinh BC=\cos(\pi/k).
$$
Hence, $\beta(f)=-4\cosh^2v-4$, where $v\in\U$, $t(v)\geq 3$.
\end{proof}

\section{Orbifolds}

Denote by $\Omega(\Gamma)$ the discontinuity set of a Kleinian group
$\Gamma$.
The {\it Kleinian orbifold}
$Q(\Gamma)=({\mathbb H}^3\cup\Omega(\Gamma))/\Gamma$
is said to be an orientable $3$-orbifold with a complete hyperbolic structure
on its interior ${\mathbb H}^3/\Gamma$ and a conformal structure
on its boundary $\Omega(\Gamma)/\Gamma$.

We need the following (Kleinian) group presentations:
\begin{itemize}
\item
$PH[\infty,m;q]=\langle x,y,s\,|\,
x^\infty=s^2=(xs)^2=(ys)^2=(xyxy^{-1})^m=(y^{-1}xys)^q=1\rangle$,
\item
$P[\infty,m;q]=\langle w,x,y,z\,|\,w^\infty=x^2=y^2=z^2=(wx)^2=(wy)^2=(yz)^2=
(zx)^q=(zw)^m=1\rangle$,
\item
$\mathcal{S}_2[\infty,m;q]=
\langle x,L\,|\,x^\infty=(xLxL^{-1})^m=(xL^2x^{-1}L^{-2})^q=1\rangle$,
\item
$GTet_1[\infty,m;q]=\langle x,y,z\,|\,
x^\infty=y^2=z^\infty=(xy)^m=(yzy^{-1}z^{-1})^q=[x,z]=1\rangle$.
\end{itemize}

Here $m$ and $q$ are integers greater than 1, or $\infty$
or $\overline\infty$ with
the following convention.
If we have a relation of the form $w^n=1$ with $n=\overline\infty$,
then we simply remove the relation $w^n=1$
from the presentation (in fact, this means that the
element $w$ is hyperbolic). Further, if $n=\infty$ and we keep
the relation $w^n=1\sim w^\infty=1$, we get a Kleinian group
presentation where parabolics are indicated.
To get an abstract group presentation, we need to remove all
relations of the form $w^\infty=1$.

\begin{theorem}\label{groups}
Let $\Gamma=\langle f,g\rangle$ be a non-elementary discrete
$\mathcal{RP}$ group
without invariant plane. Let $\beta(f)\in(-\infty,-4)$ and
let $\beta(g)=0$.
Then $\gamma(f,g)=4\cosh^2u$, where $u\in \mathcal{U}$, $t(u)\geq 3$,
and one of the following holds:
\begin{enumerate}
\item If $(t(u),2)=2$  and $\beta(f)=-4\cosh^2v/\gamma(f,g)-4$,
where $v\in \mathcal{U}$, $t(v)\geq 3$, $(t(v),2)=1$,
then $\Gamma$ is isomorphic to $PH[\infty,t(u)/2;t(v)]$.
\item If $(t(u),2)=2$ and $\beta(f)=-4\cosh^2v/\gamma(f,g)-4$,
where $v\in \mathcal{U}$, $t(v)\geq 4$, $(t(v),2)=2$,
then $\Gamma$ is isomorphic to $\mathcal{S}_2[\infty,t(u)/2;t(v)/2]$.
\item If $(t(u),2)=1$ and
$\beta(f)=-4\cosh^2v-4$, where $v\in \mathcal{U}$,
$t(v)\geq 3$, $(t(v),2)=1$,
then $\Gamma$ is isomorphic to $P[\infty,t(u);t(v)]$.
\item If $(t(u),2)=1$ and $\beta(f)=-4\cosh^2v-4$,
where $v\in \mathcal{U}$, $t(v)\geq 4$, $(t(v),2)=2$,
then $\Gamma$ is isomorphic to $GTet_1[\infty,t(u);t(v)/2]$.
\end{enumerate}
\end{theorem}

\begin{proof}
Suppose $(t(u),2)=2$, that is the dihedral angle of
$\P^*$ between $\sigma$ and $\tau$ is $2\pi/m$ with $m$ even, 
$\infty$, or $\overline\infty$.
Consider a polyhedron $\widetilde\P$ bounded by $\sigma$,
$\tau$, $\sigma'=R_\omega(\sigma)$, $\tau'=R_\omega(\tau)$,
$\eta$, and $\delta$. Applying the Poincar\'e theorem
to $\widetilde\P$ and the side pairing transformations
$g$, $g'=R_\sigma R_\omega$, $e$, and $e_f$, one can see that
$\langle g,g',e_f,e\rangle$ is isomorphic to $\widetilde\Gamma$
and has the presentation
$$
\langle f,g,e\,|\, g^\infty=e^2=(ef)^2=(eg)^2=(gfgf^{-1})^{m/2}
=(f^{-1}gfe)^p=1\rangle.
$$
If $p$ is odd, then $e\in\langle f,g\rangle$ and
$\widetilde\Gamma=\Gamma\cong PH[\infty,m/2;p]$.

If $p$ is even, $\infty$, or $\overline\infty$, 
then $\widetilde\Gamma$ contains $\Gamma$ as
a subgroup of index~$2$ and has presentation $\mathcal{S}_2[\infty,m/2;p/2]$.
In order to see this, one can apply the Poincar\'e theorem
to a polyhedron $\P$ bounded by $\tau$, $\sigma$, $\tau'$, $\sigma'$,
$\eta$, and $e_f(\eta)$, and side-pairing transformations
$f$, $g$, and $g'=fg^{-1}f^{-1}$.

The proof for $(t(u),2)=1$ is analogous. In this case we need to use
the polyhedron $\Q^*$ as the starting point.
\end{proof}

\begin{figure}[htbp]
\centering
\begin{tabular}{cc}
\includegraphics[width=2.3 cm]{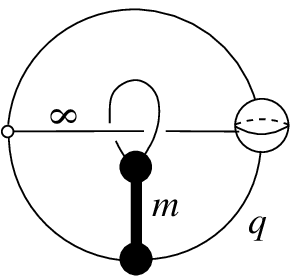} \quad & \quad
\includegraphics[width=2.3 cm]{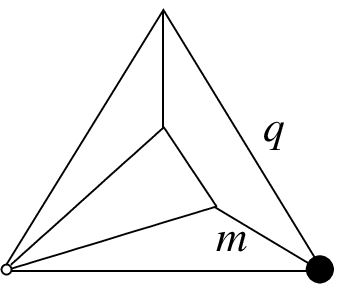}\\
(a) $\pi_1^{orb}(Q)\cong PH[\infty,m;q]$ \quad &
\quad (b) $\pi_1(Q)^{orb}\cong P[\infty,m;q]$\\
\quad $m\geq 2$, $q\geq 3$ \quad &
\quad $m\geq 3$, $q\geq 3$\\
\end{tabular}
\caption{Orbifolds embedded in ${\mathbb S}^3$}\label{ins3}
\end{figure}

\begin{figure}[htbp]
\centering
\begin{tabular}{cc}
\includegraphics[width=4.5 cm]{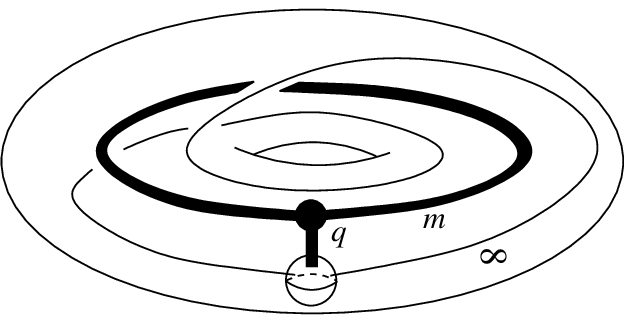} \quad & \quad
\includegraphics[width=4.5 cm]{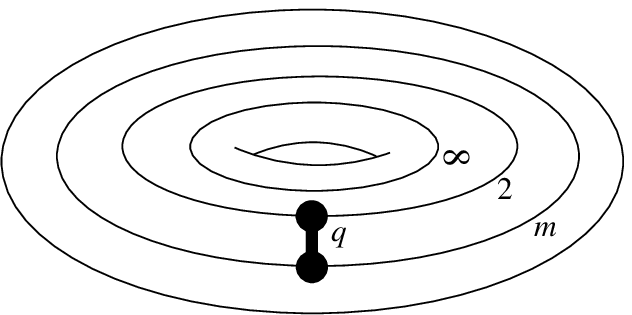}\\
(a) $\pi_1^{orb}(Q)\cong \mathcal{S}_2[\infty,m;q]$ \quad &
\quad (b) $\pi_1^{orb}(Q)\cong GTet_1[\infty,m;q]$\\
\quad $m\geq 2$, $q\geq 2$ \quad &
\quad $m\geq 3$, $q\geq 2$\\
\end{tabular}
\caption{Orbifolds embedded in Seifert fibred spaces}\label{insf}
\end{figure}

The orbifolds $Q(\Gamma)$ for the groups described in
Theorem~\ref{groups} can be obtained from corresponding fundamental
polyhedra.
In Figures~\ref{ins3} and~\ref{insf}, we schematically
draw singular sets, cusps, and boundary components of $Q(\Gamma)$
by using fat vertices and fat edges.
Roughly speaking, a fat vertex is either an interior point, or 
is removed, or removed
together with its regular neighbourhood
depending on the indices. A fat edge can be labelled by $\infty$ or
$\overline\infty$. If the index at a fat edge is $\infty$, then the egde
corresponds to a cusp, and if the index is $\overline\infty$,
the edge is removed together with its regular neighbourhood.
For details, see \cite{KK05-1}.

\medskip
In Figure~\ref{ins3}, orbifolds are embedded in ${\mathbb S}^3$ so that
$\infty$ is a non-singular interior point of $Q(\Gamma)$.
Note that the volume of $Q(PH[\infty,m;q])$ is always infinite
and $Q(P[\infty,m;q])$ is always non-compact.

\medskip
Let $T(n)$ be a Seifert fibred solid torus obtained from a trivial
fibred solid torus $D^2\times {\mathbb S}^1$ by cutting it along
$D^2\times \{x\}$ for some $x\in {\mathbb S}^1$, rotating one of the
discs through $2\pi/n$ and glueing back together.

Denote by
$\mathcal{S}(n)$ a space obtained by glueing two copies of $T(n)$
along their boundaries fibre to fibre.
Clearly, $\mathcal{S}(n)$ is homeomorphic to ${\mathbb S}^2\times {\mathbb S}^1$
and is $n$-fold covered by trivially fibred ${\mathbb S}^2\times {\mathbb S}^1$.
There are two critical fibres whose length is $n$ times shorter
than the length of a regular fibre.

\medskip

In Figure~\ref{insf}(a), orbifolds
are embedded in Seifert fibre spaces $\mathcal{S}(2)=T(2)\cup T(2)$.
We draw only the solid torus that contains
singular points (or boundary components). The other fibred torus is
meant to be attached and is not shown.
If $m<\infty$, the orbifold
$Q(\mathcal{S}_2[\infty,m;q])$
is embedded in $\mathcal{S}(2)$
in such
a manner that the axis of order $m$ 
lies on a critical fibre
of $\mathcal{S}(2)$. 
The removed regular fibre gives rise to a cusp.

In Figure~\ref{insf}(b), orbifolds
are embedded in trivially fibred space ${\mathbb S}^2\times{\mathbb S}^1$.
The rank 2 cusp corresponds to the subgroup of
$GTet_1[\infty,m;q]$ generated by $x$ and $z$.

\section{Structure of the slice $S_\infty$}

Recall that
$$
S_\infty=\{(\gamma,\beta):(\beta,0,\gamma) 
{\rm \ are\ parameters\ for\ some\ }
\langle f,g\rangle\in\D
\},
$$ 
where $\D$ denotes the class of all non-elementary discrete 
$\RP$ groups.

To investigate the slice $S_\infty$, we split the plane 
$(\gamma,\beta)$ as follows.
\begin{itemize}
\item[1.] 
If $\beta=-4$ then by \cite[Theorem~2]{KK02}, the group
$\langle f,g\rangle$ has an invariant plane. We use \cite{GGM01}
to find all discrete groups on the line $\beta=-4$.
\item[2.]
If $\beta>-4$ and $\gamma>0$ then the group $\langle f,g\rangle$
is conjugate to a subgroup of
${\rm PSL}(2,{\mathbb R})$.
More precisely, if $-4<\beta<0$ then
$f$ is elliptic and the axis of $f$ is orthogonal to an invariant
plane of $g$ and if $\beta=0$ then the fixed points of $f$ and $g$
lie in their common invariant plane. 
Discreteness criteria in terms of traces of 
$f$, $g$, and $fg$ were given in~\cite{Kna68}.
For $\beta>0$, an algorithm to decide whether
$f$ and $g$ generate a discrete group was given in \cite{GiM91}. 
\item[3.]
If $\beta>-4$ and $\gamma<0$ then $f$ is elliptic, parabolic,
or hyperbolic and the group $\langle f,g\rangle$ is known to be truly
spatial. Discrete such groups are described in \cite{KK05}, 
where $\beta$ and $\gamma$ are found explicitly.
\item[4.]
If $\beta<-4$ and $\gamma<0$ then $f$ is $\pi$-loxodromic
whose axes lies in an invariant plane of $g$. Then this plane
is invariant under action of $\langle f,g\rangle$ 
and $f$ acts as a glide-reflection on
it. A geometrical description of such discrete groups was given 
in~\cite{KS98}.
\item[5.] 
The case of $\beta<-4$ and $\gamma>0$ was treated in Section~2 of
the present paper.
\end{itemize}

We will obtain explicit formulas
for $\beta$ and $\gamma$ in the cases 2 and 4 above and 
completely describe
the structure of the slice $S_\infty$.
We will pay special attention to the subsets of $S_\infty$
corresponding to free groups.

First, we need the following elementary facts.

\begin{lemma}\label{gamma}
If $f,g\in\PSL$ and $g$ is parabolic, then
$$
\gamma(f,g)=(\tr(fg)-{\rm sign}(\tr g)\cdot\tr f)^2.
$$
\end{lemma}

\begin{proof}
By the Fricke identity, we have
\begin{eqnarray*}
\gamma(f,g)&=&\tr[f,g]-2\\
&=&\tr^2f+\tr^2g+\tr^2(fg)-\tr f\cdot \tr g\cdot \tr(fg)-4\\
&=&(\tr(fg)-{\rm sign}(\tr g)\cdot\tr f)^2,
\end{eqnarray*}
since $\tr^2g=4$.
\end{proof}

\begin{lemma}\label{fgk}
If $f,g\in\PSL$ and $\tr g=2$, then
$$
\tr(fg^k)=k(\tr(fg)-\tr f)+\tr f.
$$
\end{lemma}

\begin{proof}
By substituting $\tr g=2$ into the recurrent formula
$$\tr(fg^k)=\tr(fg^{k-1})\tr g-\tr(fg^{k-2}),$$
we immediately get the result.
\end{proof}

\begin{remark}\label{non-prim}
Suppose that $f$ is non-primitive elliptic of finite order $n$,
i.e., $\beta(f)=-4\sin^2(q\pi/n)$, where $(q,n)=1$, $1<q<n/2$.
Then there exists an integer $r$ so that $f^r$
is primitive of the same order. 
Obviously, $\langle f,g\rangle=\langle f^r,g\rangle$
and $\beta(f^r)=-4\sin^2(\pi/n)$. By \cite{GM94-1},
$\gamma(f^r,g)=(\beta(f^r)/\beta(f))\gamma(f,g)$.
\end{remark}

It is natural to introduce the constant
$$
C(q,n)=\frac{\sin^2(q\pi/n)}{\sin^2(\pi/n)}=
\frac{\beta(f)}{\beta(f^r)}\geq 1
$$
that 
plays an important role in parameters calculation concerning groups with
elliptic elements. It is also convenient to consider a parabolic
element $f$ as a limit rotation of order $n=\infty$
and write $0=\beta(f)=-4\sin^2(\pi/n)$ with $C(q,n)=C(1,n)=1$.

\subsection{$-4\leq\beta\leq 0$}

This means that $f$ is either elliptic or parabolic.
Obviously, if $f$ is elliptic of infinite order, then
$\langle f,g\rangle$ is not discrete. So we assume that
$\beta=-4\sin^2(q\pi/n)$, where $(q,n)=1$ and $1\leq q<n/2$, including
$\beta=0$.


\begin{theorem}\label{ell_rp}
Let $\Gamma=\langle f,g\rangle\subset\PSL$ have parameters
$(\beta,0,\gamma)$ with $\gamma\in{\mathbb R}\backslash \{0\}$.
Let $\beta=-4\sin^2(q\pi/n)$, where $(q,n)=1$ and $1\leq q<n/2$,
including $\beta=0$.
Then $\Gamma$ is discrete if and only if
one of the following holds:
\begin{enumerate}
\item $\gamma=-4C(q,n)\cosh^2u$, where $u\in\U$ and $t(u)\geq 3$;
\item $\gamma=4C(q,n)(\cos(\pi/n)+\cosh u)^2$, where $u\in\U$;
\item $\beta=0$ and $\gamma=4(1+\cos(2\pi/k))^2$,
where $k\geq 3$ is odd.
\end{enumerate}
\end{theorem}

\begin{proof}
Let us prove the theorem for $q=1$; in order to get the result
for $q>1$, we only need to apply Remark~\ref{non-prim}.

If $n=2$ then $\beta=-4$ and, by \cite[Theorem~4.15]{GGM01},
$\Gamma$ is discrete if and only if $\gamma=\pm 4\cosh^2u$,
where $u\in\U$ with $t(u)\geq 3$.

If $2<n\leq\infty$ and $\gamma<0$, then, by
\cite[Corollary 2.5]{KK05}, $\Gamma$ is discrete if and only if
$\gamma=-4\cosh^2u$, where $u\in\U$ and $t(u)\geq 3$.

Assume that $2<n<\infty$ and $\gamma>0$.
In this case $\Gamma$ is conjugate to a subgroup of
${\rm PSL}(2,{\mathbb R})$ and we can
apply Knapp's results \cite{Kna68} to compute $\gamma$.
Conjugate $\Gamma$ so that $\infty$ is the fixed point of $g$.
By replacing, if necessary, $f$ with $f^{-1}$ and $g$ with $g^{-1}$,
we may assume that
$$
f=\left(
\begin{array}{cc}
a& b\\
c &d
\end{array}
\right) \quad {\rm and} \quad
g=\left(
\begin{array}{rr}
-1& \tau\\
0 &-1
\end{array}
\right),
$$
where $ad-bc=1$, $a+d=-2\cos(\pi/n)$ with $n\in{\mathbb Z}$, $b>0$,
and $\tau>0$.

One can show that $\tr(fg)<2$.
By \cite[Proposition~4.1]{Kna68}, $\Gamma$ is discrete if and only if
$\tr(fg)\leq-2$ or $\tr(fg)=-2\cos(\pi/k)$, where $k\geq 2$ is an integer,
that is $\tr(fg)=-2\cosh u$, where $u\in\U$.
Hence, by Lemma~\ref{gamma},
$\gamma=(\tr(fg)+\tr f)^2=(2\cosh u+2\cos(\pi/n))^2$.

So it remains to consider the case when $n=\infty$ (i.e., $\beta=0$)
and $\gamma>0$.
Again, we normalize $\Gamma$ so that $g$ is as above and 
$f=\left(
\begin{array}{rr}
-1& 0\\
-1 &-1
\end{array}
\right)$.
By \cite[Proposition~4.2]{Kna68}, such a group
is discrete if and only if $\tau\geq 4$ or $\tau=2+2\cos(2\pi/k)$
for an integer $k\geq 3$. Since in this case $\gamma=\tau^2$, we have
that $\gamma\geq 16$ or $\gamma=(2+2\cos(2\pi/k))^2$, which can be
written as $\gamma=4(1+\cosh u)^2$, where $u\in\U$, or
$\gamma=4(1+\cos(2\pi/k))^2$ for odd $k\geq 3$.
\end{proof}

\begin{remark}
If $-4\leq\beta\leq 0$ then
$\Gamma$ is discrete and free if and only if
$\beta=0$ and $\gamma\in(-\infty,-4]\cup[16,+\infty)$.
\end{remark}

\begin{figure}[htbp]
\centering
\includegraphics[width=12.5 cm]{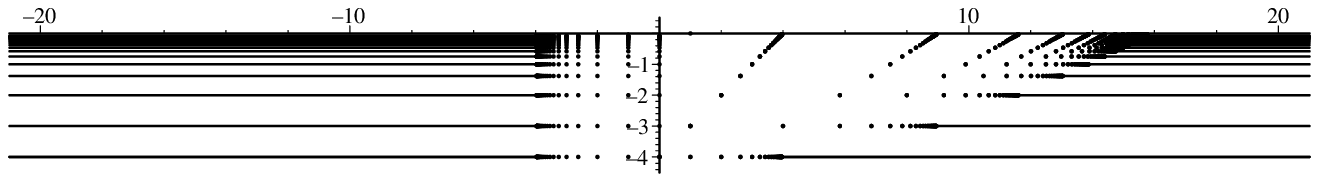}\\
$\beta=-4\sin^2(\pi/n)$ with $n\in{\mathbb Z}$\\
\bigskip
\includegraphics[width=12.5 cm]{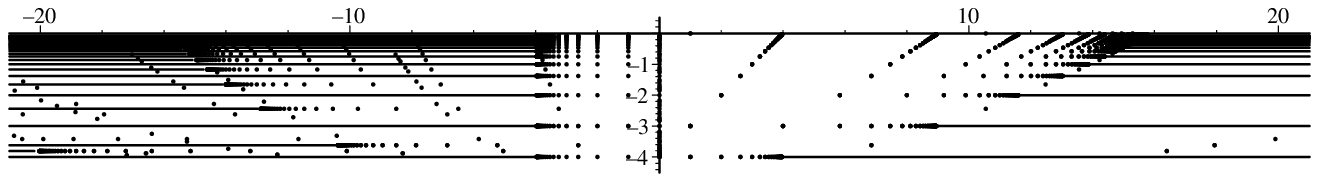}\\
$\beta=-4\sin^2(q\pi/n)$ with $(q,n)=1$ and $1\leq q<n/2$
\caption{Structure of the strip $-4\leq\beta\leq 0$}\label{strip1}
\end{figure}

The parameters from the infinite strip $-4\leq\beta\leq 0$ are
displayed in Figure~\ref{strip1}.
If $\beta=-4\sin^2(q\pi/n)$ is fixed, 
then there exist values $\gamma_1(\beta)<0$
and $\gamma_2(\beta)>0$ so that $\Gamma$ is discrete in 
the union of two rays 
$(-\infty,\gamma_1(\beta)]\cup[\gamma_2(\beta),+\infty)$.
There are only countably many discrete groups 
in $(\gamma_1(\beta),\gamma_2(\beta))$ with 
accumulation points $\gamma_1(\beta)$
and $\gamma_2(\beta)$.

Moreover, if we denote $\beta_n^q=-4\sin^2(q\pi/n)$, then
$$
\gamma_1(\beta_n^q)<\gamma_1(\beta_n^1)<\gamma_2(\beta_n^1)<
\gamma_2(\beta_n^q)\quad {\rm for\ all\ } 1<q<n/2.
$$

\subsection{$\beta>0$}

In this case $f$ is hyperbolic.

\begin{theorem}[{\cite[Corollary~2.5]{KK05}}]\label{hyp-par_rp}
Let $\Gamma=\langle f,g\rangle\subset\PSL$ have parameters
$(\beta,0,\gamma)$ with $\beta> 0$ and $\gamma<0$.
Then $\Gamma$
is discrete if and only if
$\gamma=-4\cosh^2u$, where $u\in\U$, $t(u)\geq 3$.
\end{theorem}

\begin{remark}
From \cite{KK05}, $\Gamma$ with parameters $(\beta,0,\gamma)$,
where $\beta\geq 0$ and $\gamma<0$
is free if and only if $(\gamma,\beta)$ lies in the region
$$
A=\{(\gamma,\beta):\gamma\leq -4,\beta\geq 0\}.
$$
\end{remark}

\begin{theorem}\label{hyp_fuchsian}
Let $\Gamma=\langle f,g\rangle\subset\PSL$ have parameters
$(\beta,0,\gamma)$ with $\beta>0$ and $\gamma>0$.
Let $k=\displaystyle\left\lceil
\frac{\sqrt{\beta+4}-2}{\sqrt{\gamma}}\right\rceil$.
The group $\Gamma$ is discrete if and only if one of
the following holds:
\begin{enumerate}
\item $\beta=(k\sqrt{\gamma}+2)^2-4$ and
$\gamma=16\cosh^4 u$, where $u\in\U$ and $t(u)\geq 3$;
\item $\beta=(k\sqrt{\gamma}\pm2\cos(q\pi/n))^2-4$
and
$\gamma=4C(q,n)(\cos(\pi/n)+\cosh u)^2$,
where $(q,n)=1$, $1\leq q<n/2$, and $u\in\U$;
\item $\beta=(k\sqrt{\gamma}-2\cosh u)^2-4$ and
$\gamma>4(1+\cosh u)^2$, where $u\geq 0$.
\end{enumerate}
\end{theorem}

\begin{proof}
Since $\gamma>0$, the axis of $f$ lies in
an invariant plane of $g$, so $\Gamma=\langle f,g\rangle$ is 
conjugate to a subgroup of ${\rm PSL}(2,{\mathbb R})$.
In \cite{GiM91}, an algorithm for determining whether
such a group is discrete was given.
We will apply this algorithm and calculate
parameters for each discrete group.

Normalize $\Gamma$ so that
$\infty$ is the fixed point of $g$ and $\pm1$ are the fixed points 
of~$f$. Then we can write
\begin{equation*}\label{normfg}
f=\left(
\begin{array}{cc}
a& b\\
b &a
\end{array}
\right)
\quad {\rm and}\quad
g=\left(
\begin{array}{cc}
1& \tau\\
0 &1
\end{array}
\right),
\ {\rm where}\ a^2-b^2=1,\ a>1,\ b,\tau\in{\mathbb R}.
\end{equation*}
By replacing $f$ with $f^{-1}$ and $g$ with $g^{-1}$, we may assume
that $b<0$ and $\tau>0$.

Let $k$ be a positive integer such that $\tr(fg^k)\leq 2$
and $\tr(fg^\ell)>2$ for all $\ell$ with $0\leq \ell<k$.

By Lemmas~\ref{gamma} and~\ref{fgk}, we have that
$k^2\gamma=k^2(\tr(fg)-\tr f)^2=(\tr(fg^k)-\tr f)^2$.
Since $\tr(fg^k)\leq 2$ and $\tr f>2$,
\begin{equation}\label{k2g}
\tr f=k\sqrt{\gamma}+\tr(fg^k).
\end{equation}

We distinguish three cases:

\noindent
1. $\tr(fg^k)=2$, that is $fg^k$ is parabolic.
From (\ref{k2g}),
$$
\beta=(k\sqrt{\gamma}+2)^2-4.
$$
By Theorem~\ref{ell_rp},
$\langle fg^k,g\rangle$ and, hence, $\langle f,g\rangle$
is discrete if and only if
\begin{itemize}
\item[] $\gamma=\gamma(fg^k,g)=4(1+\cosh v)^2$, where $v\in\U$, or
\item[] $\gamma=4(1+\cos(2\pi/k))^2$,
where $k\geq 3$ is odd.
\end{itemize}
These expressions can be rearranged and combined as
$\gamma=16\cosh^4u$, where $u\in\U$ and $t(u)\geq 3$.

\noindent
2. $-2<\tr(fg^k)<2$, that is $fg^k$ is elliptic and
$\tr(fg^k)=\pm 2\cos(q\pi/n)$, where $(q,n)=1$ and $1\leq q<n/2$.
Hence, from (\ref{k2g}),
$$
\beta=(k\sqrt{\gamma}\pm2\cos(q\pi/n))^2-4.
$$
By Theorem~\ref{ell_rp},
$\langle fg^k,g\rangle$ and, hence, $\langle f,g\rangle$
is discrete if and only if
$$
\gamma=4C(q,n)(\cos(\pi/n)+\cosh u)^2,
\text{\quad where } u\in\U.
$$

\noindent
3. $\tr(fg^k)\leq-2$, that is $fg^k$ is hyperbolic or parabolic
so we can write
$\tr(fg^k)=-2\cosh u$, where $u\geq 0$.
Then
$$
\beta=(k\sqrt{\gamma}-2\cosh u)^2-4.
$$
Consider the group $\langle g^{k-1}f,g\rangle$. The element
$g^{k-1}f$ is hyperbolic with $\tr(g^{k-1}f)>2$. Therefore,
one can normalize $\langle g^{k-1}f,g\rangle$ so that 
the attracting and repelling fixed points of $g^{k-1}f$ are $x_a$
and $x_r$, respectively, and $x_a<x_r$.
Since $\tr(g^kf)\leq -2$, such a group is dicrete and free 
by \cite[Case~II]{GiM91}.
So by Lemma~\ref{gamma},
we have that
\begin{eqnarray*}
\gamma=\gamma(fg^{k-1},g)&=&(\tr(fg^k)-\tr(fg^{k-1}))^2\\
&=&(2\cosh u+2\cosh v)^2,
\end{eqnarray*}
where $v$ is any positive real number.

It remains to compute $k$.
Since $\tr(fg^k)=2a+b\tau k\leq 2$, we have that
$k\geq(-2a+2)/(b\tau)$. Computing $\gamma=b^2\tau^2$, we get
$b\tau=-\sqrt{\gamma}$.
So $k=\left\lceil\displaystyle\frac{\sqrt{\beta+4}-2}
{\sqrt{\gamma}}\right\rceil$.
\end{proof}

It follows from \cite{GiM91} that $\Gamma$ is free
if and only if $(\gamma,\beta)$ lies in one of the regions
$$
C_k=\{(\gamma,\beta):\gamma\geq 16, ((k-1)\sqrt{\gamma}+2)^2\leq
\beta+4\leq (k\sqrt{\gamma}-2)^2\}, \ k=1,2,3\dots
$$

\subsection{$\beta<-4$}

First, consider $\gamma<0$. In this case the axis of the $\pi$-loxodromic
generator $f$ lies in an invariant plane of $g$ \cite{KK02}, so
$\langle f,g\rangle$ keeps this plane invariant.

\begin{theorem}\label{lox_inv}
Let $\Gamma=\langle f,g\rangle\subset\PSL$ have parameters
$(\beta,0,\gamma)$ with $\beta<-4$ and $\gamma<0$.
Let 
$k=\displaystyle\left\lceil\frac{\sqrt{-\beta-4}}{\sqrt{-\gamma}}\right\rceil$.
Then the group $\langle f,g\rangle$ is discrete if and only if
one of the following holds:
\begin{enumerate}
\item $-4(\beta+4)=\big((2k-1)\sqrt{-\gamma}\pm
\sqrt{-\gamma-8(1+\cosh u)}\big)^2$,
where $u\in\U$;
\item $4(\beta+4)=(2k-1)^2\gamma$ and $\gamma=-16\cos^2(\pi/p)$,
where $p\geq 3$ is odd;
\item $\beta=k^2\gamma-4$ and $\gamma=-4\cosh^2u$, where
$u\in\U$ and $t(u)\geq 3$.
\end{enumerate}
\end{theorem}

\begin{proof}
Let  $\delta=\{(z,t):{\rm Im}\ z=0\}$ be
the invariant plane of $\Gamma$. Since the axis of
$f$ lies in $\delta$,
we can normalize $\Gamma$ so that
the fixed point of $g$ is $\infty$, the fixed points of
$f$ are $\pm1$, and
$$
f=\left(
\begin{array}{cc}
ai& bi\\
bi &ai
\end{array}
\right),
\quad
g=\left(
\begin{array}{cc}
1& \tau\\
0 &1
\end{array}
\right),
\quad {\rm where}\ b^2-a^2=1,\ a>1,\ b,\tau\in{\mathbb R}.
$$
Further, replacing $f$ with $f^{-1}$ and $g$ with $g^{-1}$,
we can assume that $b<0$ and $\tau>0$.
Since $b$ is negative, $+1$ is the repelling fixed point of $f$ and 
$-1$ is attracting.

Let $e$ be the half-turn whose axis passes through the fixed point of $g$
orthogonally to the axis of $f$. That is $e$ fixes $0$ and $\infty$. 
Let $e_f$ and $e_1$ be half-turns such that
$f=ee_f$ and $g=e_1e$. Since $f$ is $\pi$-loxodromic, 
the axis of $e_f$ intersects the axis of $f$ (and 
the plane $\delta$) orthogonally;
denote the intersection point by $A$. 
Further, since $g$ is parabolic and 
keeps $\delta$ invariant,
the axis of $e_1$ fixes $\infty$ and lies in the plane~$\delta$.
It is easy to calculate that
$$
e=\left(
\begin{array}{rr}
i& 0\\
0 &-i
\end{array}
\right),\quad
e_f=\left(
\begin{array}{rr}
a& b\\
-b &-a
\end{array}
\right),\quad
e_1=\left(
\begin{array}{rr}
i& \tau\\
0 &-i
\end{array}
\right).
$$

Consider half-turns $e_{k-1}=g^{k-1}e$ and $e_k=g^ke$
such that $A$ lies in the region 
bounded by the axes of $e_{k-1}$
and $e_k$ in the plane $\delta$, see Figure~\ref{delta}.
It is easy to calculate that $A=-a/b-j/b$.
Since $e_k$ fixes $\infty$ and $\tau k/2$, we have that
$$
A\in\left\{(z,t):\frac{\tau(k-1)}2<{\rm Re} z\leq \frac{\tau k}2, 
\ {\rm Im}\ z=0,\ t>0\right\}.
$$ 
Hence, we can immediately determine~$k$. 
\begin{equation}\label{kbounds}
\frac{\tau(k-1)}2<-\frac ab \leq
\frac{\tau k}2.
\end{equation}
Therefore, since $2a=-i\tr f=\sqrt{-\beta-4}$ and $b\tau=-\sqrt{-\gamma}$,
$$
k=\left\lceil-\frac{2a}{b\tau}\right\rceil=
\left\lceil\frac{\sqrt{-\beta-4}}{\sqrt{-\gamma}}\right\rceil.
$$

\medskip 

It is easy to see that $\Gamma$ is discrete if and only if
$\widetilde\Gamma=\langle e_f,e_{k-1},e_k\rangle$ is. 
Following \cite{KS98},
we give geometric conditions for $\widetilde\Gamma$ to be discrete.

\medskip
\noindent
Suppose that $A\notin axis(e_k)$;
see Figure~\ref{delta}(a).
By \cite{KS98}, $\widetilde\Gamma$ is discrete if either

(a) the angle $\phi$ between $e_{k-1}$ and $e_f(e_k)$
is of the form $\pi/p$, where $p\geq 2$ is an integer, 
$\infty$, or $\overline\infty$; or

(b) $\phi=2\pi/p$, where $p\geq 3$ is odd and the bisector of
$\phi$ passes through $A$.

\smallskip
\noindent
Suppose that $A\in axis(e_k)$;
see Figure~\ref{delta}(b).
By \cite{KS98}, $\widetilde\Gamma$ is discrete if

(c) the angle $\psi$ made by $axis(e_{k-1})$ and
$axis(\tilde e_f)$ is of the form $\pi/p$, $p\geq 3$ is
an integer, $\infty$, or $\overline\infty$, where 
$\tilde e_f=e_ke_f$
is the half-turn whose axis passes through $A$ orthogonally to
$axis(e_k)$ in the plane $\delta$.

\begin{figure}[htbp]
\centering
\begin{tabular}{cc}
\includegraphics[width=5.6 cm]{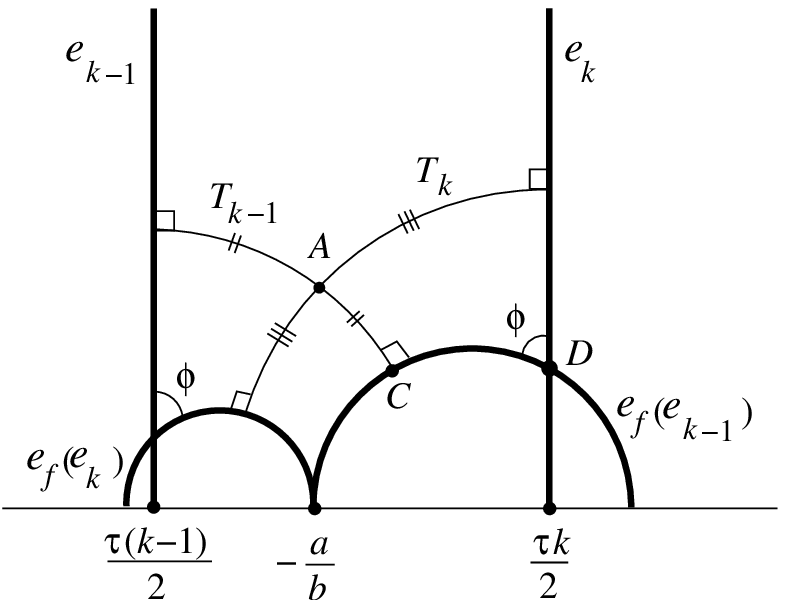}\quad
&\quad\includegraphics[width=5.2 cm]{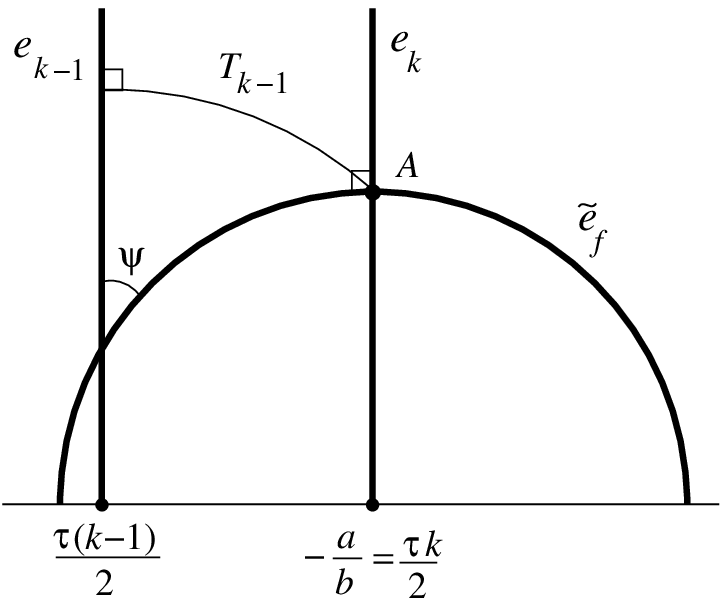}\\
(a)\quad &\quad (b)
\end{tabular}
\caption{The invariant plane $\delta$}\label{delta}
\end{figure}

There are no other discrete groups.
So, we need to calculate the parameters $\beta$ and $\gamma$
in each of the cases (a), (b), and (c).

\medskip

Assume that we are in case (a) or (b). Then each $g^\ell f=e_\ell e_f$,
$\ell\in{\mathbb Z}$,
is a $\pi$-loxodromic element with translation length $2T_\ell$
and $\tr(g^\ell f)=\pm 2i\sinh T_\ell$,
where $T_\ell$ is the distance between $e_\ell$ and~$A$.
Moreover, from the matrix representation, 
$\tr(g^\ell f)=2ai+b\tau \ell i$.
The inequalities~(\ref{kbounds})
enable us to determine the signs of $\tr(fg^{k-1})$ and $\tr(fg^k)$: 
$$
\tr(fg^k)=-2i\sinh T_k\quad {\rm and} \quad
\tr(fg^{k-1})=+2i\sinh T_{k-1}.
$$

Suppose that $p<\infty$. Simple calculations in the plane~$\delta$ 
show that
$$
\sinh CD=\frac{1+\cos\phi\cosh(2T_{k-1})}
{\sin\phi\sinh(2T_{k-1})}
$$
and, on the other hand,
$$
\sinh CD=\frac{\sinh T_k+\cos\phi\sinh T_{k-1}}
{\sin\phi\cosh T_{k-1}}.
$$
So, we obtain
$$
2(1+\cos\phi)=4\sinh T_{k-1}\sinh T_k=\tr(fg^{k-1})\tr(fg^k).
$$
Applying Lemmas~\ref{gamma} and \ref{fgk} and the facts that
$\tr f=i\sqrt{-\beta-4}$ and
$\tr(fg)-\tr f=b\tau i=-i\sqrt{-\gamma}$, we get
\begin{eqnarray*}
2(1+\cos\phi)&=&
[(k-1)(\tr(fg)-\tr f)+\tr f]\cdot[k(\tr(fg)-\tr f)+\tr f]\\
&=&k(k-1)(\tr(fg)-\tr f)^2+(2k-1)\cdot\tr f\cdot
(\tr(fg)-\tr f)+\tr^2 f\\
&=&k(k-1)\gamma+(2k-1)\sqrt{-\beta-4}\sqrt{-\gamma}+\beta+4.
\end{eqnarray*}
Hence,
$-4(\beta+4)=((2k-1)\sqrt{-\gamma}\pm\sqrt{-8(1+\cos\phi)-\gamma})^2$,
where $\phi=\pi/p$, $p\geq 2$ is an integer.
Analogous calculation can be done for $p=\infty$ and $p=\overline\infty$,
and we obtain item (1) of the theorem.

In case (b), in addition,
$T_{k-1}=T_k$. Then $\tr(fg^k)=-\tr(fg^{k-1})$ and by
Lemmas~\ref{gamma} and~\ref{fgk} we have
$$
2\sqrt{-\beta-4}=(2k-1)\sqrt{-\gamma}.
$$
Therefore,
$2(1+\cos\phi)=-\tr^2(fg^k)=(-k\sqrt{-\gamma}+\sqrt{-\beta-4})^2=
-\gamma/4$.
Hence, since $\phi=2\pi/p$, $\gamma=-16\cos^2(\pi/p)$.

Now assume that we are in case (c) and $p<\infty$.
Since in this case $e_ke_f=\tilde e_f$ is an ellitic element of order~$2$,
$\tr(g^kf)=0$. Therefore, since 
$\tr(g^kf)=-ki\sqrt{-\gamma}+i\sqrt{-\beta-4}$, we have that
$\beta=k^2\gamma-4$.

Further, since $\tr(fg^{k-1})=2i\sinh T_{k-1}$ and,
from the plane $\delta$, $\sinh T_{k-1}=\cos\psi$, we have
that
\begin{eqnarray*}
4\cos^2\psi=4\sinh^2T_{k-1}
&=&-((k-1)(\tr(fg)-\tr f)+\tr f)^2\\
&=&(-(k-1)\sqrt{-\gamma}+\sqrt{-\beta-4})^2\\
&=&(-(k-1)\sqrt{-\gamma}+k\sqrt{-\gamma})^2\\
&=&-\gamma.
\end{eqnarray*}
Thus, $\gamma=-4\cos^2(\pi/p)$,
where $p\geq 3$ is
an integer. Analogous calculations can be done for $p=\infty$
and $p=\overline\infty$ and we obtain item (3) of the theorem.
\end{proof}

\begin{remark}
If $\beta<-4$ and $\gamma<0$, then $\langle f,g\rangle$
is free if and only if $(\gamma,\beta)$ lies in one of the regions
$D_k$, $k=1,2,3,\dots$, given by
{\setlength\arraycolsep{0pt}
\begin{eqnarray*}
D_k=\{ &&(\gamma,\beta):\gamma\leq -16,\\
&&\frac{((2k-1)\sqrt{-\gamma}-\sqrt{-\gamma-16})^2}{-4}
\geq \beta+4\geq
\frac{((2k-1)\sqrt{-\gamma}+\sqrt{-\gamma-16})^2}{-4}\}.
\end{eqnarray*}
}
\end{remark}

When $\gamma>0$, the parameters were described in
Theorem~\ref{criterion_par}.
Here we just note that for $\gamma>0$ and $\beta<0$, 
the group $\langle f,g\rangle$ is free
if and only if $(\gamma,\beta)$ lies in the region
$$
B=\{(\gamma,\beta): \gamma\geq 4,\ \beta+4\leq -4/\gamma\}.
$$

\bigskip

\begin{figure}[htbp]
\begin{center}
\includegraphics[width=11 cm]{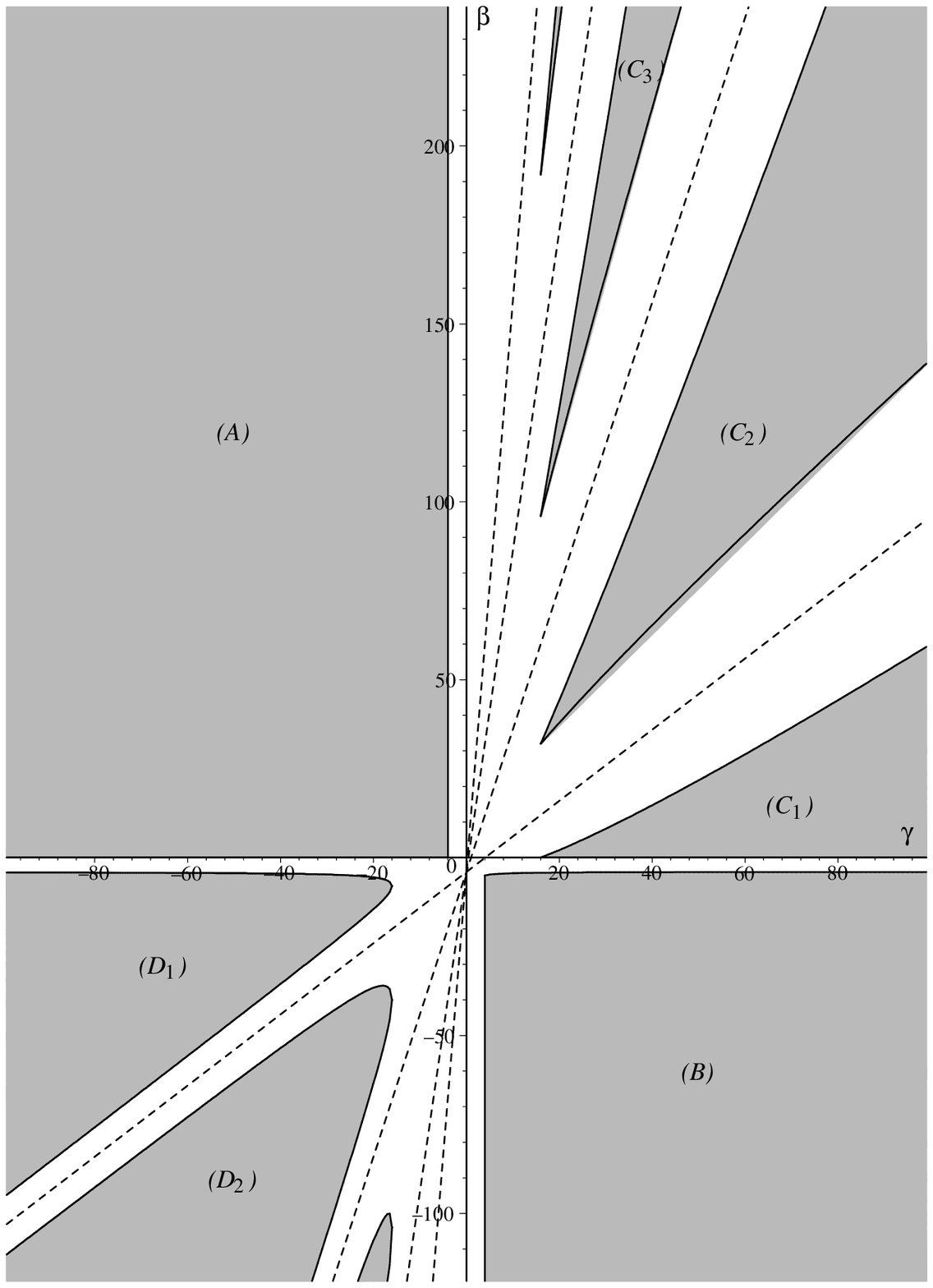}
\end{center}
{\setlength\arraycolsep{0pt}
\begin{eqnarray*}
A=\{&&(\gamma,\beta):\gamma\leq -4,\ \beta\geq 0\}\\
B=\{&&(\gamma,\beta): \gamma\geq 4,\ \beta+4\leq -4/\gamma \}\\
C_k=\{&&(\gamma,\beta):\gamma\geq 16,\ ((k-1)\sqrt{\gamma}+2)^2\leq
\beta+4\leq (k\sqrt{\gamma}-2)^2\}\\
D_k=\{&&(\gamma,\beta):\gamma\leq -16,\\
&&\frac{((2k-1)\sqrt{-\gamma}+\sqrt{-\gamma-16})^2}{-4}\leq \beta+4\leq
\frac{((2k-1)\sqrt{-\gamma}-\sqrt{-\gamma-16})^2}{-4}\}
\end{eqnarray*}
}\begin{flushleft}\quad Dashed lines $\beta=k^2\gamma-4$, 
$k=1,2,3,\dots$
\end{flushleft}
\caption{The discrete free groups}\label{free}
\end{figure}

Finally, we are able to draw those subsets of $S_\infty$
that correspond to discrete
free groups. These subsets are shown in Figure~\ref{free}.
The dashed lines $\beta=k^2\gamma-4$ are plotted to show 
a certain symmetry of $S_\infty$.

The other discrete groups contain elliptic elements.
Their parameters are represented by
lines, parabolas, hyperbolas, and points accumulating,
as orders of elliptic
elements tend to $\infty$, to the regions of
free groups.

\begin{figure}[htbp]
\centering
\includegraphics[width=11 cm]{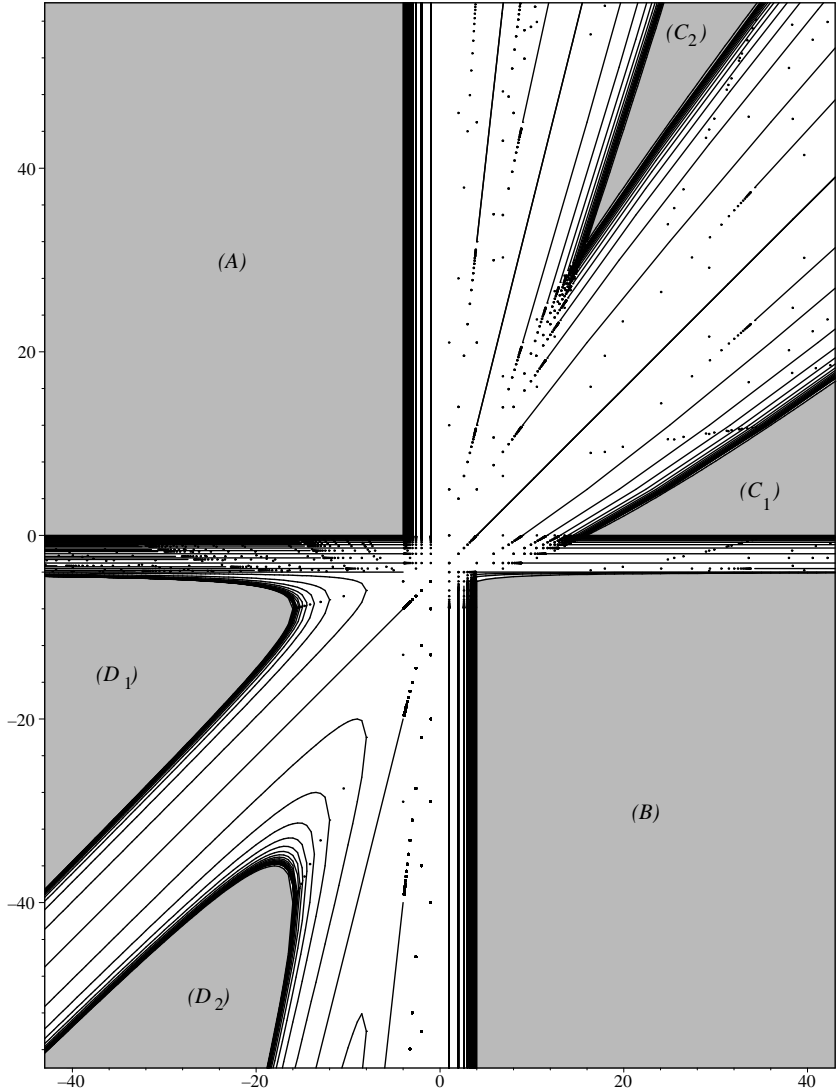}
\caption{The structure of the slice $S_\infty$}\label{whole}
\end{figure}

In Figure~\ref{whole}, the whole picture for the slice
$S_\infty$ is shown to give an idea of the structure of~$S_\infty$.
The formulas for $\beta$ and $\gamma$ obtained in Theorems
\ref{criterion_par}, \ref{ell_rp}, \ref{hyp-par_rp}, 
\ref{hyp_fuchsian}, and \ref{lox_inv},
were programmed with the package Maple~7.0
for some (sufficiently large) values of independent variables 
like $n,q\in{\mathbb Z}$ and $u,v\in\U$ and plotted on the plane 
$(\gamma,\beta)$. 

The most interesting families of parameters appear when 
$\gamma$ and $\beta$ are of the same sign.
For a fixed $k$, the hyperbolas 
$$-4(\beta+4)=\left((2k-1)\sqrt{-\gamma}\pm
\sqrt{-\gamma-8(1+\cos(\pi/p))}\right)^2,$$
where $p\geq 2$ is an integer, form a one-parameter family of curves
converging to the boundary of $D_k$ as $p\to\infty$.
Each hyperbola has the asymptotes 
$\beta=(k-1)^2\gamma-4k(1+\cos(\pi/p))+4$ and
$\beta=k^2\gamma+4k(1+\cos(\pi/p))-4$, which are obviously 
parallel to $\beta=(k-1)^2\gamma-4$ and $\beta=k^2\gamma-4$,
respectively.

\begin{figure}[htbp]
\centering
\includegraphics[width=10 cm]{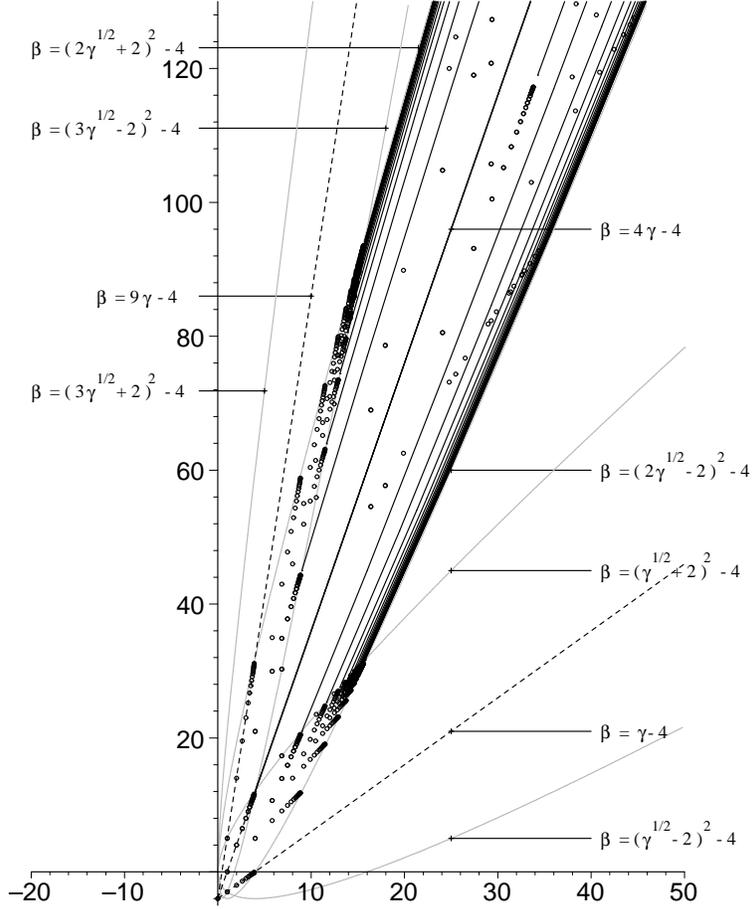}
\caption{The structure of $\Sigma_2$}\label{sigma2}
\end{figure}

For $\gamma>0$ and $\beta>0$, consider a one-parameter family of
parabolas $\beta_k=(k\sqrt{\gamma}\pm2)^2-4$. Let $\Sigma_k$
be the domain bounded by $\beta_k$:
$$
\Sigma_k=\{(\gamma,\beta): (k\sqrt{\gamma}-2)^2\leq \beta+4
\leq (k\sqrt{\gamma}+2)^2\}.
$$
Within each $\Sigma_k$, the parameters for discrete groups are given by
$$
\left\{
\begin{array}{l}
\beta=(k\sqrt{\gamma}\pm2\cos(q\pi/n))^2-4,\\
\gamma=4C(q,n)(\cos(\pi/n)+\cosh u)^2,
\end{array}
\right.
$$
where $(q,n)=1$, $1\leq q<n/2$, and $u\in\U$. 
Note that for $n=2$, we have $\beta=k^2\gamma-4$ and $\gamma=4\cosh^2u$.
As $n\to\infty$, the curves 
$\beta=(k\sqrt{\gamma}\pm2\cos(q\pi/n))^2-4$ accumulate to the
boundary of $\Sigma_k$, i.e., to the boundaries of
$C_{k-1}$ and $C_k$ (see Figure~\ref{sigma2} for an example of
$\Sigma_k$ for $k=2$).


\begin{thebibliography}{99}
%
\bibitem{Bea83}
A.~F.~Beardon, {\it The geometry of discrete groups},
Springer-Verlag, New~York--Heidelberg--Berlin, 1983.
%
\bibitem{Bea88}
A.~F.~Beardon,
{\it Fuchsian groups and $n$th roots of parabolic generators},
Holomorphic functions and moduli, Vol. II
(Berkeley, CA, 1986), 13--22,
Math. Sci. Res. Inst. Publ., 11, 1988.
%
%
\bibitem{EP94}
D.~B.~A.~Epstein and C.~Petronio, {\it An exposition of Poincar\'e's
polyhedron theorem},
L'Enseignement Math\'ematique {\bf 40} (1994), 113--170.
%
\bibitem{Fen89}
W.~Fenchel,
{\it Elementary geometry in hyperbolic space},
de Gruyter Studies in Mathematics, {\bf 11}.
Walter de Gruyter \& Co., Berlin, 1989.
%
\bibitem{GGM01}
F.~W.~Gehring, J.~P.~Gilman, and G.~J.~Martin,
{\it Kleinian groups with real parameters},
Commun. Contemp. Math. {\bf 3}, no.~2 (2001), 163--186.
%
\bibitem{GM89}
F.~W.~Gehring and G.~J.~Martin, 
{\it Stability and extremality in J\o rgensen's inequality},
Complex Variables Theory Appl. {\bf 12} (1989), no.~1-4, 277--282.
%
\bibitem{GM94-1}
F.~W.~Gehring and  G.~J.~Martin,  {\it Chebyshev polynomials and
discrete groups}, Proc. of the Conf. on Complex Analysis (Tianjin,
1992), 114--125, Conf. Proc. Lecture Notes Anal., I, Internat.
Press, Cambridge, MA, 1994.
%
\bibitem{GiM91}
J.~Gilman and B.~Maskit,  {\it An  algorithm  for  $2$-generator
Fuchsian groups}, Mich. Math.~J. {\bf 38} (1991), no.~1, 13--32.
%
\bibitem{KK02}
E.~Klimenko and N.~Kopteva, {\it Discreteness criteria for
$\mathcal{RP}$ groups},
Israel J. Math. {\bf 128} (2002), 247--265.
%
\bibitem{KK04}
E.~Klimenko and N.~Kopteva, {\it All discrete $\mathcal{RP}$ groups
whose generators have real traces},
Int. J. Algebra Comput. {\bf 15} (2005), no.~3, 577--618.
%
\bibitem{KK05}
E.~Klimenko and N.~Kopteva, {\it Discrete $\mathcal{RP}$ groups with
a parabolic generator}, 
Sib. Math. J. {\bf 46} (2005), no.~6, 1069--1076.
%
\bibitem{KK05-1}
E.~Klimenko and N.~Kopteva, {\it Two-generator Kleinian orbifolds},
2005, preprint.
%
\bibitem{KS98}
E.~Klimenko and M.~Sakuma, {\it Two-generator discrete subgroups of
${\rm Isom}({\mathbb H}^2)$
containing orientation-reversing elements}, Geometriae Dedicata {\bf
72} (1998), 247--282.
%
\bibitem{Kna68}
A.~W.~Knapp, {\it Doubly generated Fuchsian groups},  Mich. Math.
J. {\bf 15} (1968), no.~3, 289--304.
%
\bibitem{Mat82}
J.~P.~Matelski, {\it The classification of discrete 2-generator
subgroups of ${\rm PSL}(2,{\mathbb R})$}, Israel J.~Math. {\bf 42} (1982),
no.~4, 309--317.
%
\bibitem{Vin85}
E.~B.~Vinberg,
{\it Hyperbolic reflection groups},
Russian Math. Surveys {\bf 40} (1985), 31--75.
\end{thebibliography}
\end{document}